\documentclass[10pt]{elsarticle}
\journal{}
\usepackage{amsmath}
\usepackage{indentfirst}
\usepackage{amsfonts,amssymb}
\usepackage{mathrsfs}
\usepackage{amsthm}
\usepackage{lineno,hyperref}
\modulolinenumbers[5]
\newtheorem{Definition}{Definition}[section]
\newtheorem{Theorem}{Theorem}[section]

\numberwithin{equation}{section}

\begin{document}
	
	\begin{frontmatter}
		
		\title{Dynamics of an imprecise stochastic multimolecular biochemical reaction model with L\'{e}vy jumps}

		\author{Fei Sun\corref{mycorrespondingauthor}}
		\cortext[mycorrespondingauthor]{Corresponding author}
		\address{School of Mathematics and Computational Science, Wuyi University, Jiangmen 529020, China}
		\ead{sunfei@whu.edu.cn (fsun.sci@outlook.com)}

\begin{abstract}Population dynamics are often affected by sudden environmental perturbations. Parameters of stochastic models are often imprecise due to various uncertainties. In this paper, we formulate a stochastic multimolecular biochemical reaction model that includes L\'{e}vy jumps and interval parameters. Firstly, we prove the existence and uniqueness of the positive solution. Moreover, the threshold between extinction and persistence of the reaction is obtained. Finally, some simulations are carried out to demonstrate our theoretical results.
\end{abstract}


\end{frontmatter}

\section{Introduction}
\label{sec:1}

In recent decades, the study of biochemical reaction model has become one of the famous topics in mathematical biology and catalytic enzyme research. For a better review of mathematical models on the theory of biochemical reaction, see Kwek and Zhang \cite{1} and Tang and Zhang \cite{2}. 

Considering the biochemical reaction is inevitably affected by environmental noise.
Kim and Sauro \cite{3} studied the sensitivity summation theorems for stochastic biochemical reaction models.
 In order to capture essential feature of stochastic biochemical reaction systems, some researchers have
 used different methods to add random terms into the deterministic chemical reaction or epidemic models and studied the
 dynamical behavior of the corresponding stochastic models driven by white noise (see e.g. \cite{4}-\cite{11} as well as there references).

Most population systems assume that model parameters are accurately known. However, the sudden environmental perturbations may bring substantial social and economic losses. For example, the recent COVID-19 has a serious impact on the world. It is more realistic to study the population dynamics with imprecise parameters. Panja et al. \cite{23} studied a cholera epidemic model with imprecise numbers and discussed the stability condition of equilibrium points of the system. Das and Pal \cite{24} analyzed the stability of the system and solved the optimal control problem by introducing an imprecise SIR model. Other studies on imprecise parameters include those of \cite{20}-\cite{222}, and the references therein.

The main focus of this paper is Dynamics of an imprecise stochastic multimolecular biochemical reaction model with L\'{e}vy jumps. To this end, we first introduce the imprecise stochastic  multimolecular biochemical reaction model. With the help of Lyapunov functions, we prove the existence and uniqueness of the positive solution. Further, the threshold between extinction and persistence of the reaction is obtained.

The remainder of this paper is organized as follows. In Sect.~\ref{sec:2}, we introduce the basic models. In Sect.~\ref{sec:3}, the unique global positive solution of the system is proved. 
The threshold between extinction and persistence of the reaction are derived in
Sect.~\ref{sec:4} and Sect.~\ref{sec:5}. %Finally, in Sect.~\ref{sec:6}, numerical simulations and discussions are presented. \\

\section{Imprecise imprecise stochastic multimolecular biochemical reaction model}
\label{sec:2}

In this section, we introduce the imprecise stochastic multimolecular biochemical reaction model.
The multimolecular reactions described by the following reaction formulas (Selkov \cite{25}),
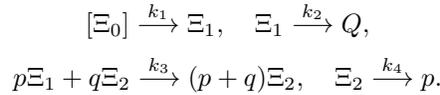
\[
[\Xi_{0} ] \stackrel{k_{1}}{\longrightarrow} \Xi_{1}, \ \ \  \Xi_{1}\stackrel{k_{2}}{\longrightarrow}  Q,
\]
\[
p\Xi_{1} + q\Xi_{2} \stackrel{k_{3}}{\longrightarrow} (p+q)\Xi_{2}, \ \ \ \Xi_{2} \stackrel{k_{4}}{\longrightarrow} p.
\]
Let $ x(t) $ and $ y(t) $ denote the concentrations of $ \Xi_{1} $ and $ \Xi_{2} $ at time $ t $, respectively, and using $ x_{0} $ to denote the concentration of $\Xi_{0}  $. Then a stochastic multimolecular biochemical reaction model with L\'{e}vy jumps takes the following form (Gao and Jiang \cite{26}).
\begin{equation}\label{2.1}
\left\{
\begin{array}{lcl}
dx = [k_{1}x_{0} - k_{2}x(t^{-}) - pk_{3}x^{p}(t^{-})y(t^{-})]dt - x^{p}(t^{-})y(t^{-})(\sigma dB(t) +  \int_{\mathbb{Y}}\gamma(u)\widetilde{N}(dt,du)),\\
dy = [pk_{3}x^{p}(t^{-})y(t^{-})- k_{4}y(t^{-})]dt + x^{p}(t^{-})y(t^{-})(\sigma dB(t) +  \int_{\mathbb{Y}}\gamma(u)\widetilde{N}(dt,du)),
\end{array}  
\right.
\end{equation}
where $ p \geq 1 $, and $ B(t) $ is standard Brownian motion with $ B(0) = 0 $. $ \sigma^{2}> 0 $ represent the intensity of white noise. 
$ \gamma(u) :\mathbb{Y} \times \Omega \rightarrow \mathbb{R}  $  is the bounded and continuous functions satisfying $ |\gamma(u)| < z $ with $ z>0 $ is a constant.
The $ x(t^{-}) $ and $ y(t^{-}) $ are the left limits of $ x(t) $ and $ y(t) $, respectively. $ \widetilde{N} $ denotes the compensated random measure defined by $ \widetilde{N}(dt, du) = N(dt, du) − \lambda(du)dt $, where
$ N $ is the Poisson counting measure and $ \lambda $ is the characteristic measure of $ N $ which is defined on a finite measurable subset $ \mathbb{Y} $ of $ (0,+\infty) $ with $ \lambda(\mathbb{Y}) < \infty $. 
We assume $ B $ and $ N $ are independent throughout the paper and denote $ \mathbb{R}^{d}_{+} = \{x \in \mathbb{R}^{d} : x_{i} > 0 \textrm{ for all } 1 \leq i \leq d\} $, $ \overline{\mathbb{R}}^{d}_{+}
= \{x \in \mathbb{R}^{d}: x_{i} \geq 0 \textrm{ for all } 1 \leq i \leq d\} $. If $ f(t) $ is an integrable function on $  [0, \infty) $, define $\langle f \rangle_{t} =\dfrac{1}{t} \int_{0}^{t} f(s) ds$.\\

Before we state the imprecise stochastic multimolecular biochemical reaction model, definitions of Interval-valued function should recalled (Pal \cite{27}).

\begin{Definition}
(Interval number) An interval number $ A $ is represented by closed interval $ [a^{l}, a^{u}] $ and defined by $ A = [a^{l}, a^{u}]  = \{x|a^{l} \leq x \leq a^{u}, x \in \mathbb{R}\},$ where $ \mathbb{R} $ is the set of real numbers and $ a^{l} $, $  a^{u} $ are the lower and upper limits of the interval numbers, respectively. The interval number $ [a, a] $ represents a real number $ a $. The arithmetic operations for any two interval numbers $ A = [a^{l}, a^{u}]$ and $ B = [b^{l}, b^{u}]$ are as follows:\\
\indent Addition: $ A+B= [a^{l}, a^{u}] + [b^{l}, b^{u}] = [a^{l}+b^{l}, a^{u}+b^{u}] $.\\
\indent Subtraction:$ A-B= [a^{l}, a^{u}] - [b^{l}, b^{u}] = [a^{l}-b^{l}, a^{u}-b^{u}] $.\\
\indent Scalar multiplication: $\alpha A= \alpha [a^{l}, a^{u}] =  [\alpha a^{l}, \alpha a^{u}]$, where $ \alpha $ is a positive real number.\\
\indent  Multiplication: $ AB= [a^{l}, a^{u}]  [b^{l}, b^{u}]= [\min \{a^{l}b^{l}, a^{u}b^{l}, a^{l}b^{u}, a^{u}b^{u} \}, \max \{a^{l}b^{l}, a^{u}b^{l}, a^{l}b^{u}, a^{u}b^{u} \}  ]  $.\\
\indent  Division: $ A/B =[a^{l}, a^{u}]/[b^{l}, b^{u}]=[a^{l}, a^{u}][ \frac{1}{b^{l}}, \frac{1}{b^{u}}] $.
\end{Definition}

\begin{Definition}
(Interval-valued function) Let $ a > 0 $, $ b > 0 $. If the interval is of the from $ [a, b] $, the interval-valued function is take as $ h(\pi) = a^{(1-\pi)}b^{\pi} $ for $ \pi \in [0, 1] $.	
\end{Definition}

Let $ \hat{k}_{1}, \hat{k}_{2}, \hat{k}_{3}, \hat{k}_{4}, \hat{p}, \hat{\sigma}$ represent the interval numbers of $ k_{1}, k_{2}, k_{3}, k_{4}, p, \sigma$, respectively. The system (\ref{2.1}) with imprecise parameters becomes,

\begin{equation}\label{2.2}
\left\{
\begin{array}{lcl}
dx = [\hat{k}_{1}x_{0} - \hat{k}_{2}x(t^{-}) - \hat{p}\hat{k}_{3}x^{\hat{p}}(t^{-})y(t^{-})]dt - x^{\hat{p}}(t^{-})y(t^{-})(\hat{\sigma} dB(t) +  \int_{\mathbb{Y}}\gamma(u)\widetilde{N}(dt,du)),\\
dy = [\hat{p}\hat{k}_{3}x^{\hat{p}}(t^{-})y(t^{-})- \hat{k}_{4}y(t^{-})]dt + x^{\hat{p}}(t^{-})y(t^{-})(\hat{\sigma} dB(t) +  \int_{\mathbb{Y}}\gamma(u)\widetilde{N}(dt,du)),
\end{array}  
\right.
\end{equation}

where $ \hat{k}_{1}= [{k_{1}}^{l}, {k_{1}}^{u}] $, $ \hat{k}_{2}= [{k_{2}}^{l}, {k_{2}}^{u}] $, $ \hat{k}_{3}= [{k_{3}}^{l}, {k_{3}}^{u}] $, $ \hat{k}_{4}= [{k_{4}}^{l}, {k_{4}}^{u}] $, $ \hat{p}= [p^{l}, p^{u}] $, $ \hat{\sigma}= [\sigma^{l}, \sigma^{u}] $.  \\

According to the Theorem 1 in Pal et al. \cite{20} and considering the interval-valued function $ f (\upsilon) = (f^{l})^{1-\upsilon}(f^{u})^{\upsilon} $ for interval $ \hat{f}= [f^{l}, f^{u}] $  for $ \upsilon \in [0, 1] $, we can prove that system (\ref{2.2}) is equivalent to the following system:

\begin{equation}\label{2.3} 
\left\{
\begin{array}{lcl}
dx = [({k_{1}}^{l})^{1-\upsilon}({k_{1}}^{u})^{\upsilon}x_{0} - ({k_{2}}^{l})^{1-\upsilon}({k_{2}}^{u})^{\upsilon}x(t^{-}) - {(p^{l})^{1-\upsilon}(p^{u})^{\upsilon}}({k_{3}}^{l})^{1-\upsilon}({k_{3}}^{u})^{\upsilon}x^{{(p^{l})^{1-\upsilon}(p^{u})^{\upsilon}}}(t^{-})y(t^{-})]dt \\
 \ \ \ \ \ \ \ - x^{{(p^{l})^{1-\upsilon}(p^{u})^{\upsilon}}}(t^{-})y(t^{-})({(\sigma^{l})^{1-\upsilon}(\sigma^{u})^{\upsilon}} dB(t) +  \int_{\mathbb{Y}}\gamma(u)\widetilde{N}(dt,du)),\\
dy = [{(p^{l})^{1-\upsilon}(p^{u})^{\upsilon}}({k_{3}}^{l})^{1-\upsilon}({k_{3}}^{u})^{\upsilon}x^{{(p^{l})^{1-\upsilon}(p^{u})^{\upsilon}}}(t^{-})y(t^{-})- ({k_{4}}^{l})^{1-\upsilon}({k_{4}}^{u})^{\upsilon}y(t^{-})]dt\\
\ \ \ \ \ \ \ + x^{{(p^{l})^{1-\upsilon}(p^{u})^{\upsilon}}}(t^{-})y(t^{-})({(\sigma^{l})^{1-\upsilon}(\sigma^{u})^{\upsilon}} dB(t) +  \int_{\mathbb{Y}}\gamma(u)\widetilde{N}(dt,du)),
\end{array}   
\right.
\end{equation}
for $ \upsilon \in [0, 1] $.
For convenience in the following investigation, let $ (\Omega, \mathcal{F}, \{\mathcal{F}_{t}\}_{t\geq 0}, \mathbb{P}) $ be a complete probability space with a filtration $ \{\mathcal{F}_{t}\}_{t\geq 0} $ satisfying the usual conditions and define
\[
R_{1} = \dfrac{{(p^{l})^{1-\upsilon}(p^{u})^{\upsilon}}({k_{3}}^{l})^{1-\upsilon}({k_{3}}^{u})^{\upsilon} }{({k_{4}}^{l})^{1-\upsilon}({k_{4}}^{u})^{\upsilon}} \Big(\dfrac{({k_{1}}^{l})^{1-\upsilon}({k_{1}}^{u})^{\upsilon}x_{0}}{({k_{2}}^{l})^{1-\upsilon}({k_{2}}^{u})^{\upsilon}}\Big)^{{(p^{l})^{1-\upsilon}(p^{u})^{\upsilon}}}
\]
and
\[
\Upsilon = \Big\{ (x,y)\in \mathbb{R}^{2}_{+}: x+y < \dfrac{({k_{1}}^{l})^{1-\upsilon}({k_{1}}^{u})^{\upsilon}x_{0}}{k} \Big\}
\]
where $ k = \min \{ ({k_{2}}^{l})^{1-\upsilon}({k_{2}}^{u})^{\upsilon}, ({k_{4}}^{l})^{1-\upsilon}({k_{4}}^{u})^{\upsilon} \} $.
We also need the following assumption.
\begin{description}
	\item[(H)] $ \Big|\big(\dfrac{({k_{1}}^{l})^{1-\upsilon}({k_{1}}^{u})^{\upsilon}x_{0}}{k}\big)^{{(p^{l})^{1-\upsilon}(p^{u})^{\upsilon}}}\gamma(u) \Big| \leq  \delta <1 \ \ \textrm{for any}\ u\in \mathbb{Y}$.
\end{description}

\begin{description}
	\item[(L)] $ R_{1}^{s} : = R_{1} - \dfrac{{\sigma''}^{2}}{2{({k_{4}}^{l})^{1-\upsilon}({k_{4}}^{u})^{\upsilon}}}(\dfrac{{({k_{1}}^{l})^{1-\upsilon}({k_{1}}^{u})^{\upsilon}}x_{0}}{k})^{2{(p^{l})^{1-\upsilon}(p^{u})^{\upsilon}}}>1 $,\\
	 where ${\sigma''}^{2} = ({{(\sigma^{l})^{1-\upsilon}(\sigma^{u})^{\upsilon}}})^{2} +  \int_{\mathbb{Y}} \dfrac{\gamma^{2}(u)}{(1-\delta)^{2}}\lambda(du)$.
\end{description}

\section{ Existence and uniqueness of positive solution of system (\ref{2.3}) }
\label{sec:3}

To research the dynamical behavior of an imprecise stochastic multimolecular biochemical reaction model, our first concern is whether the solution is global and positive. In this section, with the help of Lyapunov function, we show that system (\ref{2.3}) has a unique global positive solution with any given initial value. \\

\begin{Theorem}
	Let Assumptions (H) hold. Then for any given initial value $ (x(0), y(0)) \in \mathbb{R}^{2}_{+} $, system (\ref{2.3}) has a unique solution $ (x(t), y(t)) \in \mathbb{R}^{2}_{+} $ for all $ t \geq 0 $ almost surely (a.s.).
\end{Theorem}

\noindent \textbf{Proof.} Because the coefficients of system (\ref{2.3}) are local Lipschitz continuous (Mao \cite{28}), for any given initial value $ (x(0), y(0)) \in \mathbb{R}^{2}_{+} $, there is a unique local solution $ (x(t), y(t)) $ on $ t\in [0,\tau_{e}) $, where $ \tau_{e} $ is the explosion time (see Mao \cite{28}). In order to show that the local solution is global, we only need to prove that $ \tau_{e}=\infty $ a.s. In this context, choosing a sufficiently large number $ m_{0}\geq 1 $ such that $ (x(0), y(0)) $ lie within the interval $ [\dfrac{1}{m_{0}}, m_{0}] $. For each integer $ m \geq m_{0} $, we define the stopping time as
\[
\tau_{m}=\inf \Big\{t\in [0, \tau_{e}): \min\{x(t), y(t)\}\leq \dfrac{1}{m}   \textrm{ or }   \max \{ x(t), y(t)\}\geq m        \Big\},
\]
where $ \inf \emptyset= \infty $ ($\emptyset  $ being empty set). By the definition, $ \tau_{m} $ increases as $ m \rightarrow \infty $. Set $ \tau_{\infty}= \lim_{m \rightarrow \infty} \tau_{m} $.
Hence $\tau_{\infty}\leq \tau_{e}  $ a.s. If $ \tau_{\infty}= \infty $ a.s. is true, then $ \tau_{e}=\infty  $ a.s. for all $ t>0 $. In other words, we need to verify $ \tau_{\infty}= \infty $ a.s. If this claim is wrong, then there exist a constant $ T > 0 $ and an $ \epsilon \in (0, 1) $ such that 
\[
\mathbb{P}\{ \tau_{\infty}\leq T \} > \epsilon.
\]
Hence there is an interger $ m_{1}\geq m_{0} $ such that
\begin{equation}\label{3.1}
\mathbb{P}\{ \tau_{m}\leq T \} \geq \epsilon \textrm{ for all }  m\geq m_{1} .
\end{equation}
Consider the Lyapunov function $ V:\mathbb{R}^{2}_{+} \rightarrow \overline{\mathbb{R}}_{+} $ defined for $ (x(t), y(t)) \in \mathbb{R}^{2}_{+} $  by
\begin{equation*}
V(x,y) = (x-1-\ln x) + (y-1-\ln y),
\end{equation*}
Let $ m\geq m_{1} $ and $ T>0 $. Then, for any $ 0\leq t \leq \min{\tau_{m}, T} $, the It\v{o}’s formula (Situ \cite{29}) shows that
\begin{eqnarray*}
	dV (x, y) =& & LV(x, y)dt + {(\sigma^{l})^{1-\upsilon}(\sigma^{u})^{\upsilon}} x^{{(p^{l})^{1-\upsilon}(p^{u})^{\upsilon}}-1}y dB(t) -  {(\sigma^{l})^{1-\upsilon}(\sigma^{u})^{\upsilon}} x^{{(p^{l})^{1-\upsilon}(p^{u})^{\upsilon}}} dB(t)  \\
	&& - \int_{\mathbb{Y}}\Big[\gamma(u)x^{{(p^{l})^{1-\upsilon}(p^{u})^{\upsilon}}}y + \ln (1-\gamma(u)x^{{(p^{l})^{1-\upsilon}(p^{u})^{\upsilon}}-1}y)\Big]\widetilde{N}(dt, du)\\
	&& + \int_{\mathbb{Y}}\Big[\gamma(u)x^{{(p^{l})^{1-\upsilon}(p^{u})^{\upsilon}}}y - \ln (1+\gamma(u)x^{{(p^{l})^{1-\upsilon}(p^{u})^{\upsilon}}})\Big]\widetilde{N}(dt, du),
\end{eqnarray*}
where $ L $ is a differential operator, and
\begin{eqnarray*}
	LV(x, y) &=& (1-\dfrac{1}{x}) ({({k_{1}}^{l})^{1-\upsilon}({k_{1}}^{u})^{\upsilon}}x_{0} - {({k_{2}}^{l})^{1-\upsilon}({k_{2}}^{u})^{\upsilon}}x - {(p^{l})^{1-\upsilon}(p^{u})^{\upsilon}}{({k_{3}}^{l})^{1-\upsilon}({k_{3}}^{u})^{\upsilon}}x^{{(p^{l})^{1-\upsilon}(p^{u})^{\upsilon}}}y)\\
	&& + \dfrac{({(\sigma^{l})^{1-\upsilon}(\sigma^{u})^{\upsilon}})^{2}}{2}x^{2{(p^{l})^{1-\upsilon}(p^{u})^{\upsilon}}-2}y^{2}\\
	&&  - \int_{\mathbb{Y}}\Big[\ln(1-\gamma(u)x^{{(p^{l})^{1-\upsilon}(p^{u})^{\upsilon}}-1}y) + \gamma(u)x^{{(p^{l})^{1-\upsilon}(p^{u})^{\upsilon}}-1}y)\Big]\lambda (du)\\
	&& +  (1-\dfrac{1}{y})( {(p^{l})^{1-\upsilon}(p^{u})^{\upsilon}}{({k_{3}}^{l})^{1-\upsilon}({k_{3}}^{u})^{\upsilon}}x^{{(p^{l})^{1-\upsilon}(p^{u})^{\upsilon}}}y - {({k_{4}}^{l})^{1-\upsilon}({k_{4}}^{u})^{\upsilon}}y )\\
	&& +  \dfrac{({(\sigma^{l})^{1-\upsilon}(\sigma^{u})^{\upsilon}})^{2}}{2}x^{2{(p^{l})^{1-\upsilon}(p^{u})^{\upsilon}}} \\
	&&	- \int_{\mathbb{Y}}\Big[\ln(1+\gamma(u)x^{{(p^{l})^{1-\upsilon}(p^{u})^{\upsilon}}}) - \gamma(u)x^{{(p^{l})^{1-\upsilon}(p^{u})^{\upsilon}}}\Big]\lambda (du)\\
	&\leq&  {({k_{1}}^{l})^{1-\upsilon}({k_{1}}^{u})^{\upsilon}}x_{0} + {({k_{2}}^{l})^{1-\upsilon}({k_{2}}^{u})^{\upsilon}} + {({k_{4}}^{l})^{1-\upsilon}({k_{4}}^{u})^{\upsilon}} \\
	&&+{(p^{l})^{1-\upsilon}(p^{u})^{\upsilon}}{({k_{3}}^{l})^{1-\upsilon}({k_{3}}^{u})^{\upsilon}}x^{{(p^{l})^{1-\upsilon}(p^{u})^{\upsilon}}-1}y+ \dfrac{({(\sigma^{l})^{1-\upsilon}(\sigma^{u})^{\upsilon}})^{2}}{2}x^{2{(p^{l})^{1-\upsilon}(p^{u})^{\upsilon}}-2}y^{2} \\
	&&+  \dfrac{({(\sigma^{l})^{1-\upsilon}(\sigma^{u})^{\upsilon}})^{2}}{2}x^{2{(p^{l})^{1-\upsilon}(p^{u})^{\upsilon}}}  \\
	&& - \int_{\mathbb{Y}}\Big[\ln(1-\gamma(u)x^{{(p^{l})^{1-\upsilon}(p^{u})^{\upsilon}}-1}y) + \gamma(u)x^{{(p^{l})^{1-\upsilon}(p^{u})^{\upsilon}}-1}y)\Big]\lambda (du)\\
	&&- \int_{\mathbb{Y}}\Big[\ln(1+\gamma(u)x^{{(p^{l})^{1-\upsilon}(p^{u})^{\upsilon}}}) - \gamma(u)x^{{(p^{l})^{1-\upsilon}(p^{u})^{\upsilon}}}\Big]\lambda (du)\\
	&\leq&  {({k_{1}}^{l})^{1-\upsilon}({k_{1}}^{u})^{\upsilon}}x_{0} + {({k_{2}}^{l})^{1-\upsilon}({k_{2}}^{u})^{\upsilon}} + {({k_{4}}^{l})^{1-\upsilon}({k_{4}}^{u})^{\upsilon}} \\
	&&+ {(p^{l})^{1-\upsilon}(p^{u})^{\upsilon}}{({k_{3}}^{l})^{1-\upsilon}({k_{3}}^{u})^{\upsilon}}( \dfrac{{({k_{1}}^{l})^{1-\upsilon}({k_{1}}^{u})^{\upsilon}}x_{0}}{k} )^{{(p^{l})^{1-\upsilon}(p^{u})^{\upsilon}}}\\
	&& +  ({(\sigma^{l})^{1-\upsilon}(\sigma^{u})^{\upsilon}})^{2}( \dfrac{{({k_{1}}^{l})^{1-\upsilon}({k_{1}}^{u})^{\upsilon}}x_{0}}{k} )^{2{(p^{l})^{1-\upsilon}(p^{u})^{\upsilon}}} + H_{1} + H_{2},
\end{eqnarray*}
where 
\[
H_{1} = - \int_{\mathbb{Y}}\Big[\ln(1-\gamma(u)x^{{(p^{l})^{1-\upsilon}(p^{u})^{\upsilon}}-1}y) + \gamma(u)x^{{(p^{l})^{1-\upsilon}(p^{u})^{\upsilon}}-1}y\Big]\lambda (du),
\]
\[
H_{2} = - \int_{\mathbb{Y}}\Big[\ln(1+\gamma(u)x^{{(p^{l})^{1-\upsilon}(p^{u})^{\upsilon}}}) - \gamma(u)x^{{(p^{l})^{1-\upsilon}(p^{u})^{\upsilon}}}\Big]\lambda (du).
\]
By use of (H) and Taylor formula, we know that
\[
H_{1} = \int_{\mathbb{Y}}  \dfrac{\gamma^{2}(u)(x^{{(p^{l})^{1-\upsilon}(p^{u})^{\upsilon}}-1}y)^{2}}{2(1-\theta \gamma(u)x^{{(p^{l})^{1-\upsilon}(p^{u})^{\upsilon}}-1}y)^{2}}  \lambda (du) \leq \dfrac{\delta^{2}}{2(1-\delta)^{2}} \lambda (\mathbb{Y})
\]
with $ \theta \in (0,1) $. Similarly, we have
\[
H_{2} = \int_{\mathbb{Y}}  \dfrac{\gamma^{2}(u)(x^{{(p^{l})^{1-\upsilon}(p^{u})^{\upsilon}}})^{2}}{2(1+\theta \gamma(u)x^{{(p^{l})^{1-\upsilon}(p^{u})^{\upsilon}}})^{2}}  \lambda (du) \leq \dfrac{\delta^{2}}{2(1-\delta)^{2}} \lambda (\mathbb{Y}).
\]
Thus, we have 
\begin{equation*}
\begin{split}
LV(x, y) \leq&  {({k_{1}}^{l})^{1-\upsilon}({k_{1}}^{u})^{\upsilon}}x_{0} + {({k_{2}}^{l})^{1-\upsilon}({k_{2}}^{u})^{\upsilon}} + {({k_{4}}^{l})^{1-\upsilon}({k_{4}}^{u})^{\upsilon}}\\
&+ {(p^{l})^{1-\upsilon}(p^{u})^{\upsilon}}{({k_{3}}^{l})^{1-\upsilon}({k_{3}}^{u})^{\upsilon}}( \dfrac{{({k_{1}}^{l})^{1-\upsilon}({k_{1}}^{u})^{\upsilon}}x_{0}}{k} )^{{(p^{l})^{1-\upsilon}(p^{u})^{\upsilon}}} \\
&+ ({(\sigma^{l})^{1-\upsilon}(\sigma^{u})^{\upsilon}})^{2}( \dfrac{{({k_{1}}^{l})^{1-\upsilon}({k_{1}}^{u})^{\upsilon}}x_{0}}{k} )^{2{(p^{l})^{1-\upsilon}(p^{u})^{\upsilon}}} + \dfrac{\delta^{2}}{(1-\delta)^{2}} \lambda (\mathbb{Y}) := \widetilde{K}
\end{split}
\end{equation*}
where $\widetilde{K}  $ is a positive constant. Therefore,
\begin{equation}\label{3.2}
\begin{split}
\int_{0}^{\tau_{m} \wedge T} dV(x, y) \leq& \int_{0}^{\tau_{m} \wedge T}\widetilde{K}ds + \int_{0}^{\tau_{m} \wedge T} {{(\sigma^{l})^{1-\upsilon}(\sigma^{u})^{\upsilon}} x^{{(p^{l})^{1-\upsilon}(p^{u})^{\upsilon}}-1}(s^{-})y(s^{-})}dB(s) \\
&- \int_{0}^{\tau_{m} \wedge T} {{(\sigma^{l})^{1-\upsilon}(\sigma^{u})^{\upsilon}} x^{{(p^{l})^{1-\upsilon}(p^{u})^{\upsilon}}}(s^{-})}dB(s) \\
& - \int_{0}^{\tau_{m} \wedge T} \int_{\mathbb{Y}}\Big[\gamma(u)x^{{(p^{l})^{1-\upsilon}(p^{u})^{\upsilon}}}(s^{-})y(s^{-}) \\
&+ \ln (1-\gamma(u) x^{{(p^{l})^{1-\upsilon}(p^{u})^{\upsilon}}-1}(s^{-})y(s^{-})  )\Big] \widetilde{N}(ds, du) \\
& + \int_{0}^{\tau_{m} \wedge T} \int_{\mathbb{Y}}\Big[\gamma(u)x^{{(p^{l})^{1-\upsilon}(p^{u})^{\upsilon}}}(s^{-})y(s^{-}) - \ln (1+\gamma(u) x^{{(p^{l})^{1-\upsilon}(p^{u})^{\upsilon}}}(s^{-})  )\Big] \widetilde{N}(ds, du).
\end{split}
\end{equation}
Taking the expectations on both sides of (\ref{3.2}), we obtain that
\begin{equation}\label{3.3}
\begin{split}
\mathbb{E}V(x(\tau_{m} \wedge T ), y(\tau_{m} \wedge T) )
&\leq V(x(0), y(0) ) + \widetilde{K} \mathbb{E}(\tau_{m} \wedge T )\\
&\leq V(x(0), y(0) ) + \widetilde{K}T.
\end{split}
\end{equation}
Let $ \Omega_{m}=\{ \tau_{m} \leq T \} $ for $ m\geq m_{1} $. Then, by (\ref{3.1}), we know that $ \mathbb{P}(\Omega_{m})\geq \epsilon$. Noting that for every $ \omega\in\Omega_{m}  $, there exist $ x(\tau_{m}, \omega)$ or $ y(\tau_{m}, \omega) $, all of which equal either $ m $ or $ \frac{1}{m} $. Hence $x(\tau_{m}, \omega), y(\tau_{m}, \omega))  $ is no less than
$ m-1-\ln m $ or $ \frac{1}{m}-1-\ln \frac{1}{m} = \frac{1}{m}-1+\ln m$.
Hence
\begin{equation*}
V( x(\tau_{m}, \omega),  y(\tau_{m}, \omega))\geq (m-1-\ln m)  \wedge (\frac{1}{m}-1+\ln m).
\end{equation*}
Thus, by (\ref{3.3}), we konw that
\[
V(x(0), y(0) ) + \widetilde{K}T \geq \mathbb{E}\Big[1_{\Omega_{m}(\omega)} V(x(\tau_{m}, \omega), y(\tau_{m}, \omega)) \Big] \geq \epsilon \Big[(m-1-\ln m)  \wedge (\frac{1}{m}-1+\ln m)\Big],
\]
where $ 1_{\Omega_{m}(\omega)} $ represents the indicator function of $ \Omega_{m}(\omega) $. Setting $ m \rightarrow \infty $ leads to the contradiction
\[
\infty > V (x(0), y(0)  ) +  \widetilde{K}T = \infty.
\]
Therefore, we have  $ \tau_{\infty}= \infty $ a.s. The proof is complete. \qed

\section{Extinction}
\label{sec:4}

When studying biochemical reaction models, two of the most interesting issues are persistence and extinction. In this
section, we discuss the extinction conditions in system (\ref{2.3}) and leave its persistence to the next section.\\

\begin{Theorem}
	Let Assumption (H) hold.  For any initial value $ (x(0), y(0)) \in \Upsilon $, there is a unique positive solution $ ( x(t), y(t)) $ to system (\ref{2.3}). If one of the following two conditions holds
	\begin{description}
		\item[(A1)] $ {\sigma'}^{2} > \dfrac{({(p^{l})^{1-\upsilon}(p^{u})^{\upsilon}})^{2}{({k_{3}}^{l})^{1-\upsilon}({k_{3}}^{u})^{\upsilon}}^{2}}{2{({k_{4}}^{l})^{1-\upsilon}({k_{4}}^{u})^{\upsilon}}}$,
		\item[(A2)] $ {\sigma'}^{2} \leq \dfrac{{(p^{l})^{1-\upsilon}(p^{u})^{\upsilon}}{({k_{3}}^{l})^{1-\upsilon}({k_{3}}^{u})^{\upsilon}}}{(\dfrac{{({k_{1}}^{l})^{1-\upsilon}({k_{1}}^{u})^{\upsilon}}x_{0}}{k})^{{(p^{l})^{1-\upsilon}(p^{u})^{\upsilon}}}}$  and
		$ \dfrac{{(p^{l})^{1-\upsilon}(p^{u})^{\upsilon}}{({k_{3}}^{l})^{1-\upsilon}({k_{3}}^{u})^{\upsilon}}}{{({k_{4}}^{l})^{1-\upsilon}({k_{4}}^{u})^{\upsilon}}} (\dfrac{{({k_{1}}^{l})^{1-\upsilon}({k_{1}}^{u})^{\upsilon}}x_{0}}{k})^{{(p^{l})^{1-\upsilon}(p^{u})^{\upsilon}}}\\  -\dfrac{{\sigma'}^{2}}{2{({k_{4}}^{l})^{1-\upsilon}({k_{4}}^{u})^{\upsilon}}}(\dfrac{{({k_{1}}^{l})^{1-\upsilon}({k_{1}}^{u})^{\upsilon}}x_{0}}{k})^{2{(p^{l})^{1-\upsilon}(p^{u})^{\upsilon}}}<1 $,
	\end{description}
	then
	\[
	\limsup_{t\rightarrow \infty}\dfrac{\ln y(t)}{t} \leq   \dfrac{({(p^{l})^{1-\upsilon}(p^{u})^{\upsilon}})^{2}{({k_{3}}^{l})^{1-\upsilon}({k_{3}}^{u})^{\upsilon}}^{2}}{2{\sigma'}^{2}} - {({k_{4}}^{l})^{1-\upsilon}({k_{4}}^{u})^{\upsilon}}  < 0 \ \ a.s. \ \ \textrm{ if (A1) holds}
	\]
	and
	\[
	\limsup_{t\rightarrow \infty}\dfrac{\ln y(t)}{t} \leq  {({k_{4}}^{l})^{1-\upsilon}({k_{4}}^{u})^{\upsilon}}\Big[  \dfrac{{(p^{l})^{1-\upsilon}(p^{u})^{\upsilon}}{({k_{3}}^{l})^{1-\upsilon}({k_{3}}^{u})^{\upsilon}}}{{({k_{4}}^{l})^{1-\upsilon}({k_{4}}^{u})^{\upsilon}}} (\dfrac{{({k_{1}}^{l})^{1-\upsilon}({k_{1}}^{u})^{\upsilon}}x_{0}}{k})^{{(p^{l})^{1-\upsilon}(p^{u})^{\upsilon}}} \]
	\[
	-\dfrac{{\sigma'}^{2}}{2{({k_{4}}^{l})^{1-\upsilon}({k_{4}}^{u})^{\upsilon}}}(\dfrac{{({k_{1}}^{l})^{1-\upsilon}({k_{1}}^{u})^{\upsilon}}x_{0}}{k})^{2{(p^{l})^{1-\upsilon}(p^{u})^{\upsilon}}} -1 \Big ]  < 0 \ \ a.s. \ \ \textrm{ if (A2) holds,}
	\]
	where
	\[
	{\sigma'}^{2} =  ({(\sigma^{l})^{1-\upsilon}(\sigma^{u})^{\upsilon}})^{2} +  \int_{\mathbb{Y}} \dfrac{\gamma^{2}(u)}{(1+\delta)^{2}}\lambda(du) .
	\]
	That is to say the reaction will become extinct exponentially with probability one.
\end{Theorem}

\noindent \textbf{Proof.}
Integrating from $ 0 $ to $ t $ on both sides of (\ref{2.3}), yields
\begin{equation*}
\dfrac{x(t)-x(0)}{t} + \dfrac{y(t)-y(0)}{t} = {({k_{1}}^{l})^{1-\upsilon}({k_{1}}^{u})^{\upsilon}}x_{0} - {({k_{2}}^{l})^{1-\upsilon}({k_{2}}^{u})^{\upsilon}}\langle x \rangle_{t}- {({k_{4}}^{l})^{1-\upsilon}({k_{4}}^{u})^{\upsilon}}\langle y \rangle_{t}
\end{equation*}
Clearly, we can derive that
\begin{equation}\label{4.1}
\langle x \rangle_{t} = \dfrac{{({k_{1}}^{l})^{1-\upsilon}({k_{1}}^{u})^{\upsilon}}x_{0}}{{({k_{2}}^{l})^{1-\upsilon}({k_{2}}^{u})^{\upsilon}}} - \dfrac{{({k_{4}}^{l})^{1-\upsilon}({k_{4}}^{u})^{\upsilon}}}{{({k_{2}}^{l})^{1-\upsilon}({k_{2}}^{u})^{\upsilon}}}\langle y \rangle_{t} + \phi (t),
\end{equation}
where 
\begin{equation*}
\phi (t) =  -\dfrac{1}{{({k_{2}}^{l})^{1-\upsilon}({k_{2}}^{u})^{\upsilon}}}\Big[  \dfrac{x(t)-x(0)}{t} + \dfrac{y(t)-y(0)}{t} \Big].
\end{equation*}
This implies
\begin{equation}\label{4.2}
\lim_{t\rightarrow \infty}\phi (t) = 0 \ \ a.s. 
\end{equation}
Applying Ito formula to (\ref{2.3}) we can conclude that
\begin{equation}\label{4.3}
\begin{split}
d\ln y =& \Big[{(p^{l})^{1-\upsilon}(p^{u})^{\upsilon}}{({k_{3}}^{l})^{1-\upsilon}({k_{3}}^{u})^{\upsilon}}x^{{(p^{l})^{1-\upsilon}(p^{u})^{\upsilon}}}(t^{-})-{({k_{4}}^{l})^{1-\upsilon}({k_{4}}^{u})^{\upsilon}} - \dfrac{({(\sigma^{l})^{1-\upsilon}(\sigma^{u})^{\upsilon}})^{2}}{2} x^{2{(p^{l})^{1-\upsilon}(p^{u})^{\upsilon}}}(t^{-})\\
& + \int_{\mathbb{Y}}[\ln (1+\gamma(u)x^{{(p^{l})^{1-\upsilon}(p^{u})^{\upsilon}}}(t^{-}) ) - \gamma(u)x^{{(p^{l})^{1-\upsilon}(p^{u})^{\upsilon}}} (t^{-})] \lambda (du)\Big] dt \\
&+ {{(\sigma^{l})^{1-\upsilon}(\sigma^{u})^{\upsilon}} x^{{(p^{l})^{1-\upsilon}(p^{u})^{\upsilon}}}(t^{-})} dB(t) + \int_{\mathbb{Y}}\Big[\ln (1+\gamma(u)x^{{(p^{l})^{1-\upsilon}(p^{u})^{\upsilon}}}(t^{-}) )\Big] \widetilde{N}(dt,du).
\end{split}
\end{equation}
Integrating (\ref{4.3}) from $ 0 $ to $ t $ and then dividing by $ t $ on both sides, we obtain
\begin{equation}\label{4.4}
\begin{split}
\dfrac{\ln y(t) - \ln y(0)}{t} =& {(p^{l})^{1-\upsilon}(p^{u})^{\upsilon}}{({k_{3}}^{l})^{1-\upsilon}({k_{3}}^{u})^{\upsilon}} \langle x^{{(p^{l})^{1-\upsilon}(p^{u})^{\upsilon}}} \rangle_{t} - {({k_{4}}^{l})^{1-\upsilon}({k_{4}}^{u})^{\upsilon}} \\
&- \dfrac{({(\sigma^{l})^{1-\upsilon}(\sigma^{u})^{\upsilon}})^{2}}{2} \langle x^{2{(p^{l})^{1-\upsilon}(p^{u})^{\upsilon}}}\rangle_{t} + \dfrac{M(t)}{t} + \dfrac{\widetilde{M}(t)}{t}\\
&+ \dfrac{1}{t}\int_{0}^{t}\int_{\mathbb{Y}}\Big[\ln (1+\gamma(u)x^{{(p^{l})^{1-\upsilon}(p^{u})^{\upsilon}}}(s^{-}) ) - \gamma(u)x^{{(p^{l})^{1-\upsilon}(p^{u})^{\upsilon}}} (s^{-})\Big] \lambda (du)ds\\
\leq&  {(p^{l})^{1-\upsilon}(p^{u})^{\upsilon}}{({k_{3}}^{l})^{1-\upsilon}({k_{3}}^{u})^{\upsilon}} \langle x^{{(p^{l})^{1-\upsilon}(p^{u})^{\upsilon}}} \rangle_{t} - {({k_{4}}^{l})^{1-\upsilon}({k_{4}}^{u})^{\upsilon}}\\
& - \dfrac{1}{2}\Big[({(\sigma^{l})^{1-\upsilon}(\sigma^{u})^{\upsilon}})^{2} +  \int_{\mathbb{Y}} \dfrac{\gamma^{2}(u)}{(1+\delta)^{2}}\lambda(du)\Big]\langle x^{2{(p^{l})^{1-\upsilon}(p^{u})^{\upsilon}}}\rangle_{t}  + \dfrac{M(t)}{t} + \dfrac{\widetilde{M}(t)}{t}\\
\leq&  {(p^{l})^{1-\upsilon}(p^{u})^{\upsilon}}{({k_{3}}^{l})^{1-\upsilon}({k_{3}}^{u})^{\upsilon}} \langle x^{{(p^{l})^{1-\upsilon}(p^{u})^{\upsilon}}} \rangle_{t} - {({k_{4}}^{l})^{1-\upsilon}({k_{4}}^{u})^{\upsilon}}\\
& - \dfrac{{\sigma'}^{2}}{2} \langle x^{{(p^{l})^{1-\upsilon}(p^{u})^{\upsilon}}}\rangle_{t}^{2}  + \dfrac{M(t)}{t} + \dfrac{\widetilde{M}(t)}{t}\\
:=& h +  \dfrac{M(t)}{t} + \dfrac{\widetilde{M}(t)}{t}.
\end{split}
\end{equation}
where $ M(t): = {(\sigma^{l})^{1-\upsilon}(\sigma^{u})^{\upsilon}}\int_{0}^{t}{x^{{(p^{l})^{1-\upsilon}(p^{u})^{\upsilon}}}(s^{-})}dB(s)  $  and  $  \widetilde{M}(t):= \int_{0}^{t}\int_{\mathbb{Y}}[\ln (1+\gamma(u)x^{{(p^{l})^{1-\upsilon}(p^{u})^{\upsilon}}}(s^{-}) )]\\ \widetilde{N}(ds,du)   $ are all martingale terms and $h:(0, (\dfrac{{({k_{1}}^{l})^{1-\upsilon}({k_{1}}^{u})^{\upsilon}}x_{0}}{k})^{{(p^{l})^{1-\upsilon}(p^{u})^{\upsilon}}}) \rightarrow \mathbb{R}$ is defined by
\begin{equation}\label{4.5}
\begin{split}
h =& {(p^{l})^{1-\upsilon}(p^{u})^{\upsilon}}{({k_{3}}^{l})^{1-\upsilon}({k_{3}}^{u})^{\upsilon}} \langle x^{{(p^{l})^{1-\upsilon}(p^{u})^{\upsilon}}} \rangle_{t} - {({k_{4}}^{l})^{1-\upsilon}({k_{4}}^{u})^{\upsilon}} - \dfrac{{\sigma'}^{2}}{2} \langle x^{{(p^{l})^{1-\upsilon}(p^{u})^{\upsilon}}}\rangle_{t}^{2}\\
=&  - \dfrac{{\sigma'}^{2}}{2} \Big( \langle x^{{(p^{l})^{1-\upsilon}(p^{u})^{\upsilon}}}\rangle_{t}-\dfrac{{(p^{l})^{1-\upsilon}(p^{u})^{\upsilon}}{({k_{3}}^{l})^{1-\upsilon}({k_{3}}^{u})^{\upsilon}}}{{\sigma'}^{2}} \Big)^{2} \\
&+ \dfrac{({(p^{l})^{1-\upsilon}(p^{u})^{\upsilon}})^{2}{({k_{3}}^{l})^{1-\upsilon}({k_{3}}^{u})^{\upsilon}}^{2}}{2{\sigma'}^{2}} - {({k_{4}}^{l})^{1-\upsilon}({k_{4}}^{u})^{\upsilon}}.
\end{split}
\end{equation}
When $ {\sigma'}^{2} > \dfrac{({(p^{l})^{1-\upsilon}(p^{u})^{\upsilon}})^{2}{({k_{3}}^{l})^{1-\upsilon}({k_{3}}^{u})^{\upsilon}}^{2}}{2{({k_{4}}^{l})^{1-\upsilon}({k_{4}}^{u})^{\upsilon}}}$, (\ref{4.5}) yields 
\[
h \leq \dfrac{({(p^{l})^{1-\upsilon}(p^{u})^{\upsilon}})^{2}{({k_{3}}^{l})^{1-\upsilon}({k_{3}}^{u})^{\upsilon}}^{2}}{2{\sigma'}^{2}} - {({k_{4}}^{l})^{1-\upsilon}({k_{4}}^{u})^{\upsilon}}.
\]
Thus, by (\ref{4.4}), we have
\begin{equation}\label{4.6}
\begin{split}
\dfrac{\ln y(t)}{t} \leq&
\dfrac{\ln y(0)}{t} + h + \dfrac{M(t)}{t} + \dfrac{\widetilde{M}(t)}{t}\\
\leq& \dfrac{\ln y(0)}{t} +\dfrac{({(p^{l})^{1-\upsilon}(p^{u})^{\upsilon}})^{2}{({k_{3}}^{l})^{1-\upsilon}({k_{3}}^{u})^{\upsilon}}^{2}}{2{\sigma'}^{2}} - {({k_{4}}^{l})^{1-\upsilon}({k_{4}}^{u})^{\upsilon}} + \dfrac{M(t)}{t} + \dfrac{\widetilde{M}(t)}{t}.
\end{split}
\end{equation}
Moreover, The quadratic variation can be calculated
\[
\langle M, M \rangle_{t} = ({(\sigma^{l})^{1-\upsilon}(\sigma^{u})^{\upsilon}})^{2}\int_{0}^{t}x^{2{(p^{l})^{1-\upsilon}(p^{u})^{\upsilon}}}(s^{-})ds \leq ({(\sigma^{l})^{1-\upsilon}(\sigma^{u})^{\upsilon}})^{2} (\dfrac{{({k_{1}}^{l})^{1-\upsilon}({k_{1}}^{u})^{\upsilon}}x_{0}}{k})^{2{(p^{l})^{1-\upsilon}(p^{u})^{\upsilon}}}t,
\]
\begin{equation*}
\begin{split}
\langle \widetilde{M}, \widetilde{M} \rangle_{t} =& \int_{0}^{t} \int_{\mathbb{Y}}[\ln (1+\gamma(u)x^{{(p^{l})^{1-\upsilon}(p^{u})^{\upsilon}}}(s^{-}) )]^{2} \lambda(du)ds \\
\leq& \max \Big\{ [\ln(1-\delta)]^{2},  [\ln(1+\delta)]^{2}  \Big\} \lambda(\mathbb{Y})t.
\end{split}
\end{equation*}
Thus, by strong law of large numbers, we have
\begin{equation}\label{4.7}
\lim_{t\rightarrow \infty} \dfrac{M(t)}{t} = 0 \ \ a.s. \textrm{ and } \lim_{t\rightarrow \infty} \dfrac{\widetilde{M}(t)}{t} = 0 \ \ a.s.
\end{equation}
Clearly, taking the superior limit on both sides of (\ref{4.6}), we know that
\begin{equation*}
\limsup_{t\rightarrow \infty}\dfrac{\ln y(t)}{t} \leq \dfrac{({(p^{l})^{1-\upsilon}(p^{u})^{\upsilon}})^{2}{({k_{3}}^{l})^{1-\upsilon}({k_{3}}^{u})^{\upsilon}}^{2}}{2{\sigma'}^{2}} - {({k_{4}}^{l})^{1-\upsilon}({k_{4}}^{u})^{\upsilon}} < 0 \ \ a.s.
\end{equation*}
When $ {\sigma'}^{2} \leq \dfrac{{(p^{l})^{1-\upsilon}(p^{u})^{\upsilon}}{({k_{3}}^{l})^{1-\upsilon}({k_{3}}^{u})^{\upsilon}}}{(\dfrac{{({k_{1}}^{l})^{1-\upsilon}({k_{1}}^{u})^{\upsilon}}x_{0}}{k})^{{(p^{l})^{1-\upsilon}(p^{u})^{\upsilon}}}}$  and $ \dfrac{{(p^{l})^{1-\upsilon}(p^{u})^{\upsilon}}{({k_{3}}^{l})^{1-\upsilon}({k_{3}}^{u})^{\upsilon}}}{{({k_{4}}^{l})^{1-\upsilon}({k_{4}}^{u})^{\upsilon}}} (\dfrac{{({k_{1}}^{l})^{1-\upsilon}({k_{1}}^{u})^{\upsilon}}x_{0}}{k})^{{(p^{l})^{1-\upsilon}(p^{u})^{\upsilon}}}  -\dfrac{{\sigma'}^{2}}{2{({k_{4}}^{l})^{1-\upsilon}({k_{4}}^{u})^{\upsilon}}}(\dfrac{{({k_{1}}^{l})^{1-\upsilon}({k_{1}}^{u})^{\upsilon}}x_{0}}{k})^{2{(p^{l})^{1-\upsilon}(p^{u})^{\upsilon}}}<1 $. Then, (\ref{4.5}) yields 
\begin{equation*}
\begin{split}
h \leq& {(p^{l})^{1-\upsilon}(p^{u})^{\upsilon}}{({k_{3}}^{l})^{1-\upsilon}({k_{3}}^{u})^{\upsilon}}(\dfrac{{({k_{1}}^{l})^{1-\upsilon}({k_{1}}^{u})^{\upsilon}}x_{0}}{k})^{{(p^{l})^{1-\upsilon}(p^{u})^{\upsilon}}} \\
&- {({k_{4}}^{l})^{1-\upsilon}({k_{4}}^{u})^{\upsilon}} - \dfrac{{\sigma'}^{2}}{2}(\dfrac{{({k_{1}}^{l})^{1-\upsilon}({k_{1}}^{u})^{\upsilon}}x_{0}}{k})^{2{(p^{l})^{1-\upsilon}(p^{u})^{\upsilon}}}\\
=& {({k_{4}}^{l})^{1-\upsilon}({k_{4}}^{u})^{\upsilon}} \Big[   \dfrac{{(p^{l})^{1-\upsilon}(p^{u})^{\upsilon}}{({k_{3}}^{l})^{1-\upsilon}({k_{3}}^{u})^{\upsilon}}}{{({k_{4}}^{l})^{1-\upsilon}({k_{4}}^{u})^{\upsilon}}} (\dfrac{{({k_{1}}^{l})^{1-\upsilon}({k_{1}}^{u})^{\upsilon}}x_{0}}{k})^{{(p^{l})^{1-\upsilon}(p^{u})^{\upsilon}}}  \\
&-\dfrac{{\sigma'}^{2}}{2{({k_{4}}^{l})^{1-\upsilon}({k_{4}}^{u})^{\upsilon}}}(\dfrac{{({k_{1}}^{l})^{1-\upsilon}({k_{1}}^{u})^{\upsilon}}x_{0}}{k})^{2{(p^{l})^{1-\upsilon}(p^{u})^{\upsilon}}}-1 \Big].
\end{split}
\end{equation*}
Similarly, we get
\[
\limsup_{t\rightarrow \infty}\dfrac{\ln y(t)}{t}  \leq  {({k_{4}}^{l})^{1-\upsilon}({k_{4}}^{u})^{\upsilon}} \Big[   \dfrac{{(p^{l})^{1-\upsilon}(p^{u})^{\upsilon}}{({k_{3}}^{l})^{1-\upsilon}({k_{3}}^{u})^{\upsilon}}}{{({k_{4}}^{l})^{1-\upsilon}({k_{4}}^{u})^{\upsilon}}} (\dfrac{{({k_{1}}^{l})^{1-\upsilon}({k_{1}}^{u})^{\upsilon}}x_{0}}{k})^{{(p^{l})^{1-\upsilon}(p^{u})^{\upsilon}}}  \]
\[
-\dfrac{{\sigma'}^{2}}{2{({k_{4}}^{l})^{1-\upsilon}({k_{4}}^{u})^{\upsilon}}}(\dfrac{{({k_{1}}^{l})^{1-\upsilon}({k_{1}}^{u})^{\upsilon}}x_{0}}{k})^{2{(p^{l})^{1-\upsilon}(p^{u})^{\upsilon}}}-1 \Big] <0 \ \ a.s.
\]
which yields 
\begin{equation*}
\lim_{t\rightarrow \infty} y(t) = 0 \ \ a.s.
\end{equation*}
This completes the proof.  \qed

\section{Persistence}
\label{sec:5}

In this section, we establish sufficient conditions for persistence in the mean of system (\ref{2.3}).

\begin{Definition}\label{D51}
	Model (\ref{2.3}) is said to be persistence in the mean, if
	\[
	\liminf_{t\rightarrow \infty}\langle y \rangle_{t} >0 \ \ a.s.
	\]
\end{Definition}

\begin{Theorem}\label{T51}
	Let Assumptions (H) and (L) hold.  For any initial value $ (x(0), y(0)) \in \Upsilon $, the positive solution $ ( x(t), y(t)) $ to system (\ref{2.3}) satifies
	\begin{equation}\label{5.4}
	\begin{split}
	\liminf_{t\rightarrow \infty}\langle y \rangle_{t} \geq& -\dfrac{{({k_{2}}^{l})^{1-\upsilon}({k_{2}}^{u})^{\upsilon}}}{{({k_{4}}^{l})^{1-\upsilon}({k_{4}}^{u})^{\upsilon}}}\Big[  \dfrac{{({k_{4}}^{l})^{1-\upsilon}({k_{4}}^{u})^{\upsilon}}}{{(p^{l})^{1-\upsilon}(p^{u})^{\upsilon}}{({k_{3}}^{l})^{1-\upsilon}({k_{3}}^{u})^{\upsilon}}}\\
	& + \dfrac{{\sigma''}^{2}}{2{(p^{l})^{1-\upsilon}(p^{u})^{\upsilon}}{({k_{3}}^{l})^{1-\upsilon}({k_{3}}^{u})^{\upsilon}}}( \dfrac{{({k_{1}}^{l})^{1-\upsilon}({k_{1}}^{u})^{\upsilon}}x_{0}}{k} )^{2{(p^{l})^{1-\upsilon}(p^{u})^{\upsilon}}} \Big]^{\frac{1}{{(p^{l})^{1-\upsilon}(p^{u})^{\upsilon}}}} + \dfrac{{({k_{1}}^{l})^{1-\upsilon}({k_{1}}^{u})^{\upsilon}}x_{0}}{{({k_{4}}^{l})^{1-\upsilon}({k_{4}}^{u})^{\upsilon}}} \\
	>&0 \ \ a.s.
	\end{split}
	\end{equation}
\end{Theorem}

\noindent \textbf{Proof.}
According to (\ref{4.4}), we know that
\begin{equation}\label{5.1}
\begin{split}
\dfrac{\ln y(t) - \ln y(0)}{t} =& {(p^{l})^{1-\upsilon}(p^{u})^{\upsilon}}{({k_{3}}^{l})^{1-\upsilon}({k_{3}}^{u})^{\upsilon}} \langle x^{{(p^{l})^{1-\upsilon}(p^{u})^{\upsilon}}} \rangle_{t} - {({k_{4}}^{l})^{1-\upsilon}({k_{4}}^{u})^{\upsilon}} \\
&- \dfrac{({(\sigma^{l})^{1-\upsilon}(\sigma^{u})^{\upsilon}})^{2}}{2} \langle x^{2{(p^{l})^{1-\upsilon}(p^{u})^{\upsilon}}}\rangle_{t} + \dfrac{M(t)}{t} + \dfrac{\widetilde{M}(t)}{t}\\
&+ \dfrac{1}{t}\int_{0}^{t}\int_{\mathbb{Y}}\Big[\ln (1+\gamma(u)x^{{(p^{l})^{1-\upsilon}(p^{u})^{\upsilon}}}(s^{-}) ) - \gamma(u)x^{{(p^{l})^{1-\upsilon}(p^{u})^{\upsilon}}} (s^{-})\Big] \lambda (du)ds\\
\geq&  {(p^{l})^{1-\upsilon}(p^{u})^{\upsilon}}{({k_{3}}^{l})^{1-\upsilon}({k_{3}}^{u})^{\upsilon}} \langle x^{{(p^{l})^{1-\upsilon}(p^{u})^{\upsilon}}} \rangle_{t} - {({k_{4}}^{l})^{1-\upsilon}({k_{4}}^{u})^{\upsilon}}\\
& - \dfrac{({(\sigma^{l})^{1-\upsilon}(\sigma^{u})^{\upsilon}})^{2}}{2}( \dfrac{{({k_{1}}^{l})^{1-\upsilon}({k_{1}}^{u})^{\upsilon}}x_{0}}{k} )^{2{(p^{l})^{1-\upsilon}(p^{u})^{\upsilon}}}\\
&-( \dfrac{{({k_{1}}^{l})^{1-\upsilon}({k_{1}}^{u})^{\upsilon}}x_{0}}{k} )^{2{(p^{l})^{1-\upsilon}(p^{u})^{\upsilon}}} \int_{\mathbb{Y}} \dfrac{\gamma^{2}(u)}{2(1-\delta)^{2}}\lambda(du)
+ \dfrac{M(t)}{t} + \dfrac{\widetilde{M}(t)}{t}\\
=&  {(p^{l})^{1-\upsilon}(p^{u})^{\upsilon}}{({k_{3}}^{l})^{1-\upsilon}({k_{3}}^{u})^{\upsilon}} \langle x^{{(p^{l})^{1-\upsilon}(p^{u})^{\upsilon}}} \rangle_{t} - {({k_{4}}^{l})^{1-\upsilon}({k_{4}}^{u})^{\upsilon}} 
-\dfrac{1}{2}\Big[({(\sigma^{l})^{1-\upsilon}(\sigma^{u})^{\upsilon}})^{2}\\
& +  \int_{\mathbb{Y}} \dfrac{\gamma^{2}(u)}{(1-\delta)^{2}}\lambda(du)\Big]( \dfrac{{({k_{1}}^{l})^{1-\upsilon}({k_{1}}^{u})^{\upsilon}}x_{0}}{k} )^{2{(p^{l})^{1-\upsilon}(p^{u})^{\upsilon}}}  + \dfrac{M(t)}{t} + \dfrac{\widetilde{M}(t)}{t}\\
=& {(p^{l})^{1-\upsilon}(p^{u})^{\upsilon}}{({k_{3}}^{l})^{1-\upsilon}({k_{3}}^{u})^{\upsilon}} \langle x^{{(p^{l})^{1-\upsilon}(p^{u})^{\upsilon}}} \rangle_{t} - {({k_{4}}^{l})^{1-\upsilon}({k_{4}}^{u})^{\upsilon}}\\
& - \dfrac{1}{2}{\sigma^{''}}^{2}( \dfrac{{({k_{1}}^{l})^{1-\upsilon}({k_{1}}^{u})^{\upsilon}}x_{0}}{k} )^{2{(p^{l})^{1-\upsilon}(p^{u})^{\upsilon}}}  + \dfrac{M(t)}{t} + \dfrac{\widetilde{M}(t)}{t}.
\end{split}
\end{equation}
This together with (\ref{4.1}) implies
\begin{equation}\label{5.2}
\begin{split}
\dfrac{\ln y(t) - \ln y(0)}{t} \geq& {(p^{l})^{1-\upsilon}(p^{u})^{\upsilon}}{({k_{3}}^{l})^{1-\upsilon}({k_{3}}^{u})^{\upsilon}} \Big[ \dfrac{{({k_{1}}^{l})^{1-\upsilon}({k_{1}}^{u})^{\upsilon}}x_{0}}{{({k_{2}}^{l})^{1-\upsilon}({k_{2}}^{u})^{\upsilon}}}- \dfrac{{({k_{4}}^{l})^{1-\upsilon}({k_{4}}^{u})^{\upsilon}}}{{({k_{2}}^{l})^{1-\upsilon}({k_{2}}^{u})^{\upsilon}}}\langle y \rangle_{t} + \phi(t)  \Big]^{{(p^{l})^{1-\upsilon}(p^{u})^{\upsilon}}}\\
&-  {({k_{4}}^{l})^{1-\upsilon}({k_{4}}^{u})^{\upsilon}} - \dfrac{1}{2}{\sigma^{''}}^{2}( \dfrac{{({k_{1}}^{l})^{1-\upsilon}({k_{1}}^{u})^{\upsilon}}x_{0}}{k} )^{2{(p^{l})^{1-\upsilon}(p^{u})^{\upsilon}}}  + \dfrac{M(t)}{t} + \dfrac{\widetilde{M}(t)}{t}.
\end{split}
\end{equation}
Since $ -\infty < \ln y(t) < \ln \dfrac{{({k_{1}}^{l})^{1-\upsilon}({k_{1}}^{u})^{\upsilon}}x_{0}}{k} $, (\ref{5.2}) yields
\begin{equation}\label{5.3}
\begin{split}
\langle y \rangle_{t} \geq& -\dfrac{{({k_{2}}^{l})^{1-\upsilon}({k_{2}}^{u})^{\upsilon}}}{{({k_{4}}^{l})^{1-\upsilon}({k_{4}}^{u})^{\upsilon}}}\Big\{   \dfrac{1}{{(p^{l})^{1-\upsilon}(p^{u})^{\upsilon}}{({k_{3}}^{l})^{1-\upsilon}({k_{3}}^{u})^{\upsilon}}}\big[   {({k_{4}}^{l})^{1-\upsilon}({k_{4}}^{u})^{\upsilon}} \\
&+ \dfrac{1}{2}{\sigma^{''}}^{2}( \dfrac{{({k_{1}}^{l})^{1-\upsilon}({k_{1}}^{u})^{\upsilon}}x_{0}}{k} )^{2{(p^{l})^{1-\upsilon}(p^{u})^{\upsilon}}} +\dfrac{\ln \frac{{({k_{1}}^{l})^{1-\upsilon}({k_{1}}^{u})^{\upsilon}}x_{0}}{k} - \ln y(0)}{t} -\dfrac{M(t)}{t} - \dfrac{\widetilde{M}(t)}{t} \big]  \Big\}^{\frac{1}{{(p^{l})^{1-\upsilon}(p^{u})^{\upsilon}}}}\\
& + \dfrac{{({k_{1}}^{l})^{1-\upsilon}({k_{1}}^{u})^{\upsilon}}x_{0}}{{({k_{4}}^{l})^{1-\upsilon}({k_{4}}^{u})^{\upsilon}}} + \dfrac{{({k_{2}}^{l})^{1-\upsilon}({k_{2}}^{u})^{\upsilon}}}{{({k_{4}}^{l})^{1-\upsilon}({k_{4}}^{u})^{\upsilon}}}\phi(t).
\end{split}
\end{equation}
Taking the inferior limit on both sides of (\ref{5.3}) and combining with
Lemma~, from (\ref{4.2}) and (\ref{4.7}) we have
\begin{equation*}
\begin{split}
\liminf_{t\rightarrow \infty}\langle y \rangle_{t} \geq& -\dfrac{{({k_{2}}^{l})^{1-\upsilon}({k_{2}}^{u})^{\upsilon}}}{{({k_{4}}^{l})^{1-\upsilon}({k_{4}}^{u})^{\upsilon}}} \Big[  \dfrac{{({k_{4}}^{l})^{1-\upsilon}({k_{4}}^{u})^{\upsilon}}}{{(p^{l})^{1-\upsilon}(p^{u})^{\upsilon}}{({k_{3}}^{l})^{1-\upsilon}({k_{3}}^{u})^{\upsilon}}} \\
&+  \dfrac{1}{2{(p^{l})^{1-\upsilon}(p^{u})^{\upsilon}}{({k_{3}}^{l})^{1-\upsilon}({k_{3}}^{u})^{\upsilon}}}{\sigma^{''}}^{2}( \dfrac{{({k_{1}}^{l})^{1-\upsilon}({k_{1}}^{u})^{\upsilon}}x_{0}}{k} )^{2{(p^{l})^{1-\upsilon}(p^{u})^{\upsilon}}}  \Big]^{\frac{1}{{(p^{l})^{1-\upsilon}(p^{u})^{\upsilon}}}}\\
& + \dfrac{{({k_{1}}^{l})^{1-\upsilon}({k_{1}}^{u})^{\upsilon}}x_{0}}{{({k_{4}}^{l})^{1-\upsilon}({k_{4}}^{u})^{\upsilon}}}
\ \ a.s.
\end{split}
\end{equation*}
Therefore, by Assumption (L), we can easily obtain (\ref{5.4}). \qed

\end{document}